\documentclass{article}
\usepackage{amsmath}
\usepackage{amssymb}
\usepackage{amsthm}
\usepackage{latexsym}
\usepackage{epsfig}
\usepackage{url}

\begin{document}
\newcommand{\ontop}[2]{\genfrac{}{}{0pt}{}{#1}{#2}}

\newtheorem{theorem}{Theorem}[section]
\newtheorem{lemma}[theorem]{Lemma}
\newtheorem{proposition}[theorem]{Proposition}
\newtheorem{corollary}[theorem]{Corollary}
\newtheorem{definition}[theorem]{Definition}

\title{Selecting universities: personal preference and rankings}

\author{Peter Huggins, Lior Pachter \\ Department of
Mathematics\\ University of California at Berkeley}

\date{\today}
\maketitle

{\abstract{
Polyhedral geometry can be used to quantitatively assess the dependence of rankings on personal preference, 
and provides a tool for both students and universities to assess US News and World Report rankings.     
}}

\section{Introduction}

The annual US News and World Report (USNWR) rankings have become a major factor influencing college admissions \cite{cite1}.  
University placements affect the decisions of students, and the rankings are carefully monitored by the institutions that are compared.  Administrators tout the rankings when they are useful for recruitment and fundraising \cite{cite2} while at the same time lamenting perceived flaws in the methods used to gather data and establish rankings. They also bemoan the fluctuating pressures exerted on their institutions as rankings change from year to year \cite{cite3}.   

The current USNWR rankings of top national universities are based on aggregate scores computed for 130 universities.  The scores consist of linear combinations of measurements of several attributes, such as prestige, selectivity, student graduation and retention, faculty resources, and student satisfaction \cite{cite4}.  Criticisms of USNWR have addressed many of the ingredients that contribute to the final rankings.  For example, there has been widespread disagreement on what type of data should be collected.  Furthermore, even when there has been some agreement that certain attributes of universities are informative for comparison, there has been disagreement on how to perform the measurements.  Finally, there is still no consensus on how to combine the various scores to obtain a single number reflecting the overall appeal of a university.   

This last point is the focus of our study.  We maintain that any information collected about a university is valuable, although its relative importance to individuals may differ based on their personal preferences.  Certain students may value the prestige afforded to selective schools, while others may be concerned with student/faculty ratios or class size.  Similarly, universities may wish to advertise their strengths to targeted groups who may best appreciate them.  The need to aggregate scores is based on the difficulty in comparing multiple features for many schools.  However instead of a single weighting, we believe it is useful and informative to examine multiple options. Indeed, USNWR has admitted as much.  In a National Opinion Research Center study funded by USNWR to evaluate its methodology \cite{cite5} it was reported that "The principal weakness of the current approach is that the weights used to combine the various measures into an overall rating lack any defensible empirical or theoretical basis."   

According to USNWR methodology \cite{cite4}, seven main categories of attributes comprise the final aggregate score: 25\% peer assessment, 20\% freshman retention, 20\% faculty resources, 15\% student selectivity, 10\% financial resources, 5\% graduation performance, and 5\% alumni giving. Various measurements for each attribute form the basis for assigning each national university a numerical score in each of the seven categories.  The seven numerical scores are then weighted and summed with the above listed weights to give the university's final score.  We studied how universities' final scores and rankings depend on the relative weights assigned to the seven categories.  We restricted our attention to nonnegative weights only, since each of the seven categories is measuring a desirable feature. We note that our methods are not restricted to just seven categories, although computational limits may prohibit analyses for too many attributes. In such cases, the methodology of \cite{cite6} can be used to first identify correlates among attributes and thereby reduce the dimension of the problem.  

Our study required the seven category scores for each university.  USNWR publishes three of the category scores, and reports only the rankings for the other four categories \cite{cite7}. We therefore reverse-engineered the four unknown categories, from the final scores and category rankings published by USNWR. This was accomplished by casting the problem in terms of feasibility for a linear program. We also simulated random data for the seven categories and repeated our analysis to verify that our findings are typical for rankings of multidimensional data. 

The mathematical problem of determining the sensitivity of rankings to changes in the weights is the geometric problem of determining the vertices and normal fans of certain polytopes called $k$-ordertopes.  A seven-dimensional vector describes the seven attribute scores assigned to each university.  The convex hull of all universities' 7d vectors is a polytope, whose vertices are all the possible top-ranked universities that can arise for some choice of weights.  Similarly, the $k$-ordertope is a polytope that can be constructed whose vertices represent the top $k$ universities that can arise for some choice of weights.The normal fan of a $k$-ordertope gives the complete partition of scoring weights according to the top $k$ ranking they induce.

We computed the non-negative envelopes of $k$-ordertopes for all 7 categories in USNWR for $k \leq r$.  For $k=4$ this computation involved finding the non-negative envelope of more than 270 million points in seven-dimensional space. This was accomplished using Theorem 2.7 and Lemma 4.1, in addition to novel software. We found that there are a total of 27 universities that can appear among the top four for some choice of nonnegative weights.  However, only nine appear in the top four for more than 5\% of the possible weights. These are Princeton, Harvard, Yale, Stanford, University of Pennsylvania, Caltech, MIT, Duke, and Dartmouth.  Any generic ranking is therefore likely to identify these top schools. In \cite{cite1, cite8} it was reported that rankings influence college choices especially for top-ranked universities. Thus, it is noteworthy that the order of the top-ranked universities can vary for a significant fraction of weights. For example, we found that Yale (USNWR \#3) beats Princeton (USNWR \#1) for 26\% of weight choices. Although exact order is sensitive to weight choices, both Yale and Princeton (and Harvard) are ranked in the top 5 for 95 percent of weight choices; the USNWR data confirms they are top schools irrespective of the weights chosen. However for many universities, no single ranking is informative.  We computed each university's ranking interval (RI), which is the smallest range of ranks for the university which are induced by at least 95\% of the possible choices of scoring weights. The average interval width of a university's RI was 35.5; for universities in USNWR's top 25 the average RI width was 11.1, and the average RI width for the remaining 105 universities was 41.3.  Some schools, such as Georgia Tech (RI of [25, 73]), have extremely large RIs compared to similarly ranked schools, meaning that their specific placement by USNWR is essentially arbitrary.  

The ranking interval measures the smallest range of ranks induced by at least 95\% of possible scoring weights.   As the average RI interval widths show, significant instability in ranks is the rule rather than the exception. We found that all but the top three universities ranked by USNWR (Princeton, Harvard, Yale) had an RI of six or more ranks.  This confirms that while the top universities are good all around regardless of the scoring weights, the rankings of other universities depend strongly on individual preferences.   

The proliferation of college rankings \cite{cite9} has created a confusing situation where multiple rankings that are difficult to compare and contrast are circulated, frequently without the underlying raw data. $k$-ordertopes organize attribute score vectors in a way that reveals the dependence of rankings on weights. They allow for analysis of the robustness of individual placements, and for aggregate statistics to be compiled. Furthermore, the examination of an entire $k$-ordertope provides a novel way to compare and contrast universities. The volumes of the normal fans of $k$-ordertopes also provide a quantitative basis for selecting weights that match desired rankings. 

We believe that the data collected by USNWR is useful for both students and universities, but ideally USNWR will publish the raw measurements and leave the rankings to the reader, as is done for German universities by the CHE/DAAD consortium \cite{cite10}.  The CHE/DAAD consortium provides a website where users can explore the strengths and weaknesses of universities by specifying preferences in weights before rankings are displayed online.  Furthermore, in the USNWR magazine a few different rankings could be published for a handful of popular preferences. The publication of a single ranking provides readers with a summary that may be misleading and fails to properly account for personal preference and uncertainty in weights.  By examining $k$-ordertopes, students can also account for the uncertainty in their preferences in a systematic way, and universities can identify the spectrum of weights where they rank highly so that they can better target prospective students.

\section{Polyhedral geometry}

In this section we introduce the basic definitions and theorems from polyhedral geometry necessary for our results. We refer the reader to \cite{Ziegler} for more details.

   {\definition{
The {\bf affine hull} of a set of points $V = \{v_1, \ldots, v_n\} \subset {\bf R}^d$ is the set $\{\sum_{i=1}^n c_i v_i \,\, | \,\, \sum_{i=1}^n c_i = 1\}$. 
}}

{\definition{
The {\bf convex hull} of a set of points $V = \{v_1, \ldots, v_n\} \subset {\bf R}^d$ is the set $\{\sum_{i=1}^n c_i v_i \,\, | \,\, \sum_{i=1}^n c_i = 1, \,\, \forall c_i \geq 0\}$. 
}}

{\definition
A {\bf polytope} is a convex hull any finite non-empty $V \subset {\bf R}^d$.
}

Equivalently, a polytope can be defined as a bounded polyhedron, where a polyhedron is the set of solutions of finitely many weak linear inequalities.  

The {\em dimension} of a polytope $P \subset {\bf R^d}$ is the dimension of its affine hull as a manifold.  To avoid confusion between $d$ and $\dim P$, $d$ is called the {\em ambient dimension} of $P$.

{\definition
Given a polytope $P \subset {\bf R}^d$ and a functional $c \in {\bf R}^d$, the {\bf face} $F_c \subset P$ is the set $\hbox{\em argmax}_{x \in P} c \cdot x$.
}

Any face of a polytope is again a polytope.  Note that $P$ is a face of itself by putting $c = 0$.  For generic choices of $c$, the face $F_c$ will be a single point, which is called a {\em vertex} of $P$.  1-dimensional faces are called {\em edges}, $(\dim P - 1)$-dimensional faces are called {\em facets}, and 
$(\dim P - 2)$-dimensional faces are called {\em ridges}.  
%The set of all faces of $P$ form a finite distributive lattice (ordered by inclusion), which is called the {\em face lattice} of $P$.

{\definition
Given a polytope $P \subset {\bf R}^d$ and a face $F \subset P$, the {\bf normal cone} of $F$ is the set of all functionals
$c$ for which $F_c = F$.
}

The relative interiors of normal cones of all faces of $P$ partition ${\bf R}^d$, and form a lattice ordered by inclusion.  
The set of normal cones has other special properties as well, and is called the {\em normal fan} of $P$.

{\definition
Given a polytope $P \subset {\bf R}^d$ and a matrix $A \in {\bf R}^{m \times d}$, 
let $C = \{c \,\,\, | \,\,\, Ac \geq 0 \}$.  The {\bf $A$-envelope} of $P$ is the set of
faces $\{F \subset P \, | \, N(F) \cap C \neq 0 \}$.  If $A = I$, then the $A$-envelope is called the {\bf non-negative
envelope} of $P$.
}

The main objects we compute and analyze are polytopes called {\em $k$-ordertopes}. For each $k$, the vertices of the $k$-ordertope correspond to top $k$ university rankings which are attainable for some
choice of weights. The $k$-ordertope is defined as follows:  Consider a set of vectors $u_i$ that are distinct (we work with the score vectors from the {\em US News} data).
  Fix positive decreasing weights $\alpha_1 > \alpha_2 \ldots > \alpha_n > 0$, where $n$ is the number of universities.
  The {\em $k$-ordertope} $P_k$ is the convex hull of all weighted sums of $k$ universities' score vectors:
  \[
    P_k = \hbox{ conv}\{ \sum_{i = 1}^k \alpha_i u_{\sigma_i}\}_{\sigma_1, \ldots, \sigma_k}          
  \]
  where $\sigma_1, \sigma_2, \ldots \sigma_k$ range over all possible $n(n-1)\cdot(n-k+1)$ choices of $k$ distinct indices.
  Each vertex of the $k$-ordertope corresponds to a choice of top $k$
  universities which is obtainable for at least one choice of scoring scheme $c$:

  {
  \begin{theorem} 
   A point $\sum_{i = 1}^k \alpha_i u_{\sigma_i} \in P_k$ is a vertex of $P_k$ if and only if 
   there is linear scoring scheme $c$ which satisfies $c \cdot u_{\sigma_1} > c \cdot u_{\sigma_2} > \ldots > c \cdot u_{\sigma_k}$, 
   and $c \cdot u_{\sigma_k} > c \cdot u_j$ for all $j \notin \{\sigma_1, \ldots, \sigma_k\}$.
    The normal fan of $P_k$ is independent of the choice of weights $\alpha_i$.   
  \end{theorem}
  }
 
  {\em Proof.}  Suppose $c$ satisfies  $c \cdot u_{\sigma_1} > \ldots > c \cdot u_{\sigma_k}$, 
   and $c \cdot u_{\sigma_k} > c \cdot u_j$ for all $j \notin \{\sigma_1, \ldots, \sigma_k\}$.  
Then we claim that 
   $c \cdot \sum_{i = 1}^k \alpha_i u_{\sigma_i} > c \cdot \sum_{i = 1}^k \alpha_i u_{\gamma_i}$, unless $\gamma_i = \sigma_i$ for all
   $i$.  By basic convexity, this would imply that $\sum_{i = 1}^k \alpha_i u_{\sigma_i}$ is a vertex on $P_k$.

   So consider any choice of distinct $\gamma_1, \gamma_2, \ldots \gamma_k$.

   \begin{itemize}
  
   \item If $c \cdot u_{\gamma_i} < c \cdot u_{\gamma_j}$ for some $i < j$, then since $\alpha_i > \alpha_j$ we have 
  \[
     c \cdot (\alpha_i u_{\gamma_j} + \alpha_j u_{\gamma_i}) > c \cdot (\alpha_i u_{\gamma_i} + \alpha_j u_{\gamma_j})
  \]

   So by swapping the values of $\gamma_i$ and $\gamma_j$, we can obtain a new point in $P_k$ which has a higher dot-product
   $c \cdot \sum_{i = 1}^k \alpha_i u_{\gamma_i}$.

   \item If $c \cdot u_{\gamma_1} > c \cdot u_{\gamma_2}  > \ldots > c \cdot u_{\gamma_k}$, and if
         $\gamma_k \notin \{\sigma_1, \ldots \sigma_k\}$, then there must also be some 
         $\sigma_i \notin \{\gamma_1, \ldots, \gamma_k\}$ by basic counting principles.  
         We have $c \cdot u_{\sigma_i} \geq c \cdot u_{\sigma_k} > c \cdot u_{\gamma_k}$, and hence, by redefining
         $\gamma_k := \sigma_i$, we can obtain a new point in $P_k$ which increases the dot-product 
         $c \cdot \sum_{i = 1}^k \alpha_i u_{\gamma_i}$.

   \end{itemize}

   Thus, $\sum_{i = 1}^k \alpha_i u_{\gamma_i}$ cannot maximize the functional $c$ over $P_k$, unless 
   $c \cdot u_{\gamma_1} > c \cdot u_{\gamma_2}  > \ldots > c \cdot u_{\gamma_k}$ and 
   $\gamma_k \in \{\sigma_i\}$.  But then {\em all} $\gamma_i$ must be contained in $\{\sigma_i\}$,
   and so $\{\gamma_i\} = \{\sigma_i\}$.  Furthermore, since $c \cdot u_{\gamma_1} > \ldots > c \cdot u_{\gamma_k}$ and 
   $c \cdot u_{\sigma_1} > \ldots > c \cdot u_{\sigma_k}$, we must also have $\gamma_i = \sigma_i$ for all $i$.

   This completes the proof of the `if' part of the theorem.  The `only if' follows from the `if':  If $\sum_{i = 1}^k \alpha_i u_{\sigma_i} \in P_k$ is a vertex, then it must strictly maximize some functional $c$, and this $c$ must induce an ordering 
  $c \cdot u_{\gamma_1} \geq \ldots \geq c \cdot u_{\gamma_n}$ on the score vectors, where the first $k$ inequalities are strict.   Then $\sum_{i = 1}^k \alpha_i u_{\gamma_i}$  is the vertex which maximizes the functional $c$ 
   over $P_k$, via the `if' part of the theorem.  The proof of the `if' part also shows that $\sum_{i = 1}^k \alpha_i u_{\sigma_i}$
   cannot equal this vertex unless $\sigma_i = \gamma_i$ for all $i=1,\ldots,k$, as desired.  

   In fact we have proved that $c$ is in the relative interior of a vertex normal cone $N(\sum_{i=1}^k \alpha_i u_{\sigma_i})$ iff $c \cdot u_{\sigma_1} > \ldots > c \cdot u_{\sigma_k}$ 
   and $c \cdot u_{\sigma_k} > c \cdot u_j$ for all $j \notin \{\sigma_1, \ldots, \sigma_k\}$.  As a corollary, we have that the interiors of vertices' normal cones depend only on the $u_i$ and are independent of the choice of $\alpha_i$.  This implies that the entire normal fan is independent of the choice of $\alpha_i$ as well.  \qed

  Our theorem characterizes normal cones of $P_k$ and shows why they are important.  For each choice of $k$ distinct
  universities $u_{\sigma_1}, \ldots, u_{\sigma_k}$, the normal cone $N(\sum_{i=1}^{k} \alpha_i u_{\sigma_i})$ gives all possible functionals
  which make $u_{\sigma_i}$ the $i$th ranked university for $i = 1,2,\ldots k$.  Thus $P_k$ and its normal fan completely describe how the top $k$ rankings
  depend on the choice of functional $c$.

  We are only interested in {\em reasonable} vertices, i.e. vertices whose normal cones' interiors intersect the 
  non-negative first orthant ${\bf R}^7_{\geq 0}$.  These vertices give precisely the top $k$ rankings which are attainable
  under a reasonable scoring scheme, in which no categories are penalized.  Since we are not interested in any other vertices, we compute the non-negative 
  envelope of $k$-ordertopes, instead of the entire $k$-ordertopes.

  \section{Data}

  We studied the {\em US News and World Report}'s Best Colleges 2008 rankings of national 
universities.  This data set comprises numerical scores 
 for each of 130 universities, in seven categories:

  \begin{itemize}
   \item Peer assessment (university officials are asked to rate the 130 universities on a scale from 1 to 5),
   \item Freshman retention rate,
   \item Graduation performance (number of percentage points above or below the university's expected graduation rate as predicted by USNWR),
    \item Faculty resources (class size, faculty to student ratio, faculty salaries),
    \item Student selectivity (SAT scores, application rejection rate),
    \item Financial resources,
     \item Alumni giving (percentage of alumni that make donations).
   \end{itemize}
The data can therefore be described by 7-dimensional vectors of measurements ($u_i$, $i = 1, \ldots, 130$), which we call the {\em score vectors.} In some cases individual category scores are obtained by combining several different measurements. We do not analyze the data at that level of detail, and focus on the $7$ categories that are used to produce the overall scores.

  In order to aggregate a university's score vector into an overall score for the university, 
  a linear scoring scheme is used
  where each university $u_i$ is assigned a score of $c \cdot u_i$, for some fixed choice of $c \in {\bf R}^7$.  
  All categories under consideration  are assumed to be intrinsically `good', 
  so $c$ is assumed to have non-negative entries.  We can also assume, without loss of generality, that the entries
  in $c$ sum to 1.  US News uses the specific choice of $c = (0.25, 0.20, 0.05, 0.20, 0.15, 0.10, 0.05)$.

 USNWR does not make available the score vectors $u_i$, but it does release three of the components (peer assessment, graduation performance, and alumni giving), as well as the universities' overall scores $c \cdot u_i$.  USNWR also publishes
the rankings of the universities in each of the four unpublished categories.

  \section{Methods}

  \subsection{Algorithm to compute non-negative envelopes of $k$-ordertopes}

  Our definition of the $k$-ordertope $P_k$ provides an oracle to find a vertex of $P_k$ which maximizes a given functional $c$.  Namely, the oracle sorts
  all $u_i$ according to the value of $c  \cdot u_i$, and then returns $\sum_{i = 1}^k \alpha_i u_i$.  (In case of ties $c \cdot u_i = c \cdot u_j$, the lexicographically
  smaller of $u_i$ and $u_j$ is ranked higher.)

  Given such an oracle for finding vertices, the software {\tt iB4e} \cite{iB4e} builds the entire polytope by carefully choosing functionals and 
  querying the oracle to find vertices, until all vertices are guaranteed to be found.  {\tt iB4e} constructs both the vertices and facets of the final polytope.  It runs in polynomial time when the dimension of the polytope is fixed, and the number of oracle queries is exactly $V + F$, where $V$ is the number
 of vertices of the constructed polytope and $F$ is the number of facets.
  
  Since we are only interested in reasonable vertices, we modify the oracle so that it can also return points of the form $-N e_j$, where $e_j$ is a standard
  basis vector in ${\bf R}^7$, and $N$ is a large positive number.  Then the outputted polytope is the convex hull of $P_k \cup \{-Ne_1, \ldots, -Ne_7\}$ instead of $P_k$.  
  If $N$ is sufficiently large, then this polytope gives the non-negative envelope of $P_k$, as the following lemma shows.

  {\lemma
  Let $P \subset {\bf R}^d$ be a $d$-dimensional polytope, and let $\{ e_j \}$ be the standard basis.  Define $Q_N = conv(P \cup \{ -Ne_1, \ldots -Ne_d \} )$.  If $N$ is sufficiently large, and if $v \in P$ is a vertex such that $N(v) \cap {\bf R}^d_{>0} \neq 0$, then $v$  is a vertex of $Q_N$. If $N$ is sufficiently large and $v \in P$ is a vertex of $Q_N$, then $v$ is a vertex on the non-negative envelope of $P$.
  }
 
  {\em Proof.}  Suppose $v$ is a vertex of $P$ for which there is $c \in N(v) \cap {\bf R}^d_{>0}$.  
Let $\epsilon = \min_i c_i$.  Then $\epsilon > 0$ and we have $c(-Ne_j) < -N \epsilon$ for all $-Ne_j$.  So if $N$ is large enough, we have $cv > c(-Ne_j)$ for all $-Ne_j$, and thus $v$ is 
  a vertex of not only $P$ but also $Q_N$.  This proves the first half of the theorem.
  
  Now suppose $v \in P$ is not a vertex on the non-negative envelope of $P$.  Then $v$ cannot weakly maximize any non-zero functional $c \in {{\bf R}_{\geq 0}}^d$.  By the Farkas Lemma, there must be a point $x \in P$ such that $x_j > v_j$ for all $j = 1, \ldots d$.  Thus $v = x + \sum_j \epsilon_j(-Ne_j)$, where all $\epsilon_j$ are positive.  We rewrite this as an affine combination of $x$ and the $-Ne_j$:

  \begin{equation}\label{eq:1}
     v = (1 - \delta)x + \sum_j (\epsilon_j - \delta x_j / N) (-Ne_j) 
  \end{equation}

  where $\delta = (\sum_j \epsilon_j) / (1 + \sum_j (x_j/N))$.  If $N$ is large, then the $\epsilon_j$ are small, and $\delta \simeq \sum_j \epsilon_j$ will also be small (and positive since the $\epsilon_j$ are positive).  Furthermore if $N$ is large then all coefficients in equation (\ref{eq:1}) for $v$ are positive, 
so equation (\ref{eq:1}) expresses $v$ as a strictly convex combination of $x$ and the $-Ne_j$.  Thus 
  $v$ is not a vertex of $Q_N$ for large $N$ (in fact, $v$ is an interior point).  $\square$
 
  \subsection{Measuring normal cones}     
  
  Absent any canonical choice of scoring scheme $c$, we measure the fraction of non-negative scoring schemes which give each ranking.  
  This means measuring normal cones of $k$-ordertopes, restricted to the first orthant. Since all functionals along a common ray 
  give the same ranking, we measure a cone with respect to the Gaussian density $f$ restricted to the first orthant:  

  \[
    f(c) \sim e^{{-||c||}^2} {\hbox{ if }} c > 0, \,\,\, f(c) = 0 {\hbox{ otherwise }}   
  \]

  Geometrically, this is 
  equivalent to intersecting a cone with the unit hypersphere and first orthant, and then computing the volume of the intersection region, which is a {\em spherical polytope}.

  Computing exact volumes of spherical polytopes is difficult in dimensions four or higher.  But we can easily sample from the first orthant with respect to the Gaussian density.  
  This allows us to approximate cone volumes quite well for all normal cones whose volume is not too small.  
  This is sufficient for our purposes since we are not interested in accurately approximating the probabilities of highly unlikely rankings.
  We can in fact accurately approximate the volume of every normal cone though, by approximating a spherical polytope by a simplicial mesh.  
  
  This gives a rational alternative to using just one (subjective) choice of $c$.  For example, for two universities $u_1, u_2$ in the $1$-ordertope,
  if $N(u_1) \cap {\bf R_{\geq 0}}^7$ is larger than $N(u_2) \cap {\bf R_{\geq 0}}^7$, then there is a larger percentage of choices of $c$ which make 
  $u_1$ the top-ranked university instead of $u_2$.

  We use volumes of normal cones to compute ranking intervals for each university.  The 95\% ranking interval for university $u_i$ is the shortest contiguous interval of ranks $I$ such that at least 95\% of functionals assign $u_i$ a ranking in $I$.  In case a university has multiple 95\% ranking intervals of the same width, we choose the leftmost interval containing the smallest ranks.
 
  \subsection{Reverse-engineering of unpublished categories}

  For each university, the seven category scores computed by USNWR give a vector $u_i \in {\bf R}^7$.  We only know $u_{i1}, u_{i3}$ and 
  $u_{i7}$ corresponding to peer assessment, graduation performance, and alumni giving.  However USNWR does publish rankings in the other four categories.
  For example Harvard (USNWR overall rank \#2) is ranked \#1 in retention and Princeton (USNWR rank \#1) is ranked \#2 in retention.  Thus, even though we don't know the retention scores $u_{22}$ for Harvard and $u_{12}$ for Princeton, we do know that Harvard's retention score is greater than Princeton's, i.e. $u_{22} \geq u_{12} + \epsilon$ for some fixed $\epsilon > 0$.  USNWR also publishes the overall final score $c \cdot u_i$ for each university.

  Thus, there is a system of linear equations and inequalities that the unknown category scores must satisfy.  The set of solutions to these linear constraints
  is a polytope.  However, since there are 130 universities
  and four unknown categories, the polytope of possible values for the unknown scores is several-hundred dimensional, so we cannot compute the entire polytope explicitly. However standard linear programming can quickly find vertices of the polytope that maximize given functionals.  

  We could use linear programming to find a vertex of the polytope and use the vertex as an estimate of the unknown category scores.  However, vertices are 
  extremal points of the polytope, whereas ideally we want interior points. In theory there are fast approximate sampling algorithms for high-dimensional polytopes, but they are complicated and to our knowledge none have been implemented.  Instead we computed many vertices and averaged them
  to obtain a central point.

In order to assess the quality of our predicted category scores, we treated one of the three known categories (alumni giving) as an unknown and therefore were able to use it as a control in order to measure the quality of our estimates. We therefore estimated category scores for five unknown categories.
  We minimized and maximized 
  each university's unknown category score separately, which required solving $5 \times 130 \times 2 = 1300$ linear programs.  We then averaged the 1300 vertices to obtain our estimate for the five unknown categories.  

  We then compared our estimate $\hat{a} \in {\bf R}^{130}$ for alumni giving to the true alumni giving scores $a$.  After standardizing so that $a, \hat{a}$ have
  standard deviation 1, we found that $|a - \hat{a}|$ was $0.0991$ on average. This is not surprising because $a$ and $\hat{a}$ are in the same order --- so they are already highly correlated.  To see if our estimate $\hat{a}$ was significantly good, we sampled random values for the entries in $\hat{a}$ from the normal distribution with the same mean and standard deviation as $a$, and then sorted $\hat{a}$ to be in the same order as $a$.  Out of 10000 trials, 95\% of the random $\hat{a}$ produced a value of $mean|a - \hat{a}| > 0.168$, the average value of $mean|a - \hat{a}|$ was $0.236$, and the lowest value of $mean|a - \hat{a}|$ was $0.10210$.

  \section{Computational results}

  Using our reverse-engineered estimates of the USNWR data, we computed the non-negative envelope of $k$-ordertopes for $k$ up to 4.  For $k = 4$, our computation entails finding the non-negative envelope
  of $130 \cdot 129 \cdot 128 \cdot 127 > 270,000,000$ points.  Such a large computation would be essentially impossible by conventional means of listing the points and determining which points lie on
  the envelope.  Thus our use of vertex-finding oracles and {\tt iB4e} was cruical to accomplish the computation.   The number of vertices of the $k$-ordertopes which we found on the boundary of the non-negative envelope are given in Table \ref{table:verts}.  
  \begin{table}
  \label{table:verts}
\begin{center}
  \begin{tabular}{ll}
   $k$ & \# Vertices \\
\hline
    1 &  10 \\ 
    2 &  59 \\
    3 & 276 \\ 
    4 & 1082 \\
   \end{tabular}
\end{center}
\caption{Number of vertices of the $k$-ordertopes.}
\end{table}
  Note that there are over a thousand possible top 4 rankings that USNWR could have conceivably published by using a different
  choice of non-negative weights.

 The schools which can be ranked \#1 are Princeton, Harvard, Yale, Stanford, U Penn, Cal-tech, MIT, Notre Dame, Johns Hopkins, and Penn State.  All of these are USNWR top 20 schools except Penn State (USNWR rank \#48).  Penn State is able to be ranked \#1 by virtue of its high graduation performance---which is a category weighted very low in USNWR's scoring scheme.

  We further found that 17 schools can appear in the top 2, 26 schools can appear in the top 3, and 27 schools can appear in the top 4.  

We show an example of part of the normal fan of a 3-ordertope in Figure 1.  In the example we have restricted the analysis to three of the features weighted highest by USNWR (peer assessment, retention, student selectivity) so that the picture is 3-dimensional.  Each choice of weights is a vector that lies inside one of the cones that form the normal fan of the 3-ordertope. The top 3 universities for each choice of weights label the cones. For example, with a high preference for selectivity, the ranking will be 1. Harvard 2. Yale 3. Princeton, whereas a combined preference for peer assessment and freshman retention results in 1. Harvard 2. Princeton 3. Stanford. The normalized volumes of the cones measure the fraction of weights for which each ranking is optimal.

\begin{figure}[ht]
\label{fig:kset_korder}
\includegraphics[scale=0.65]{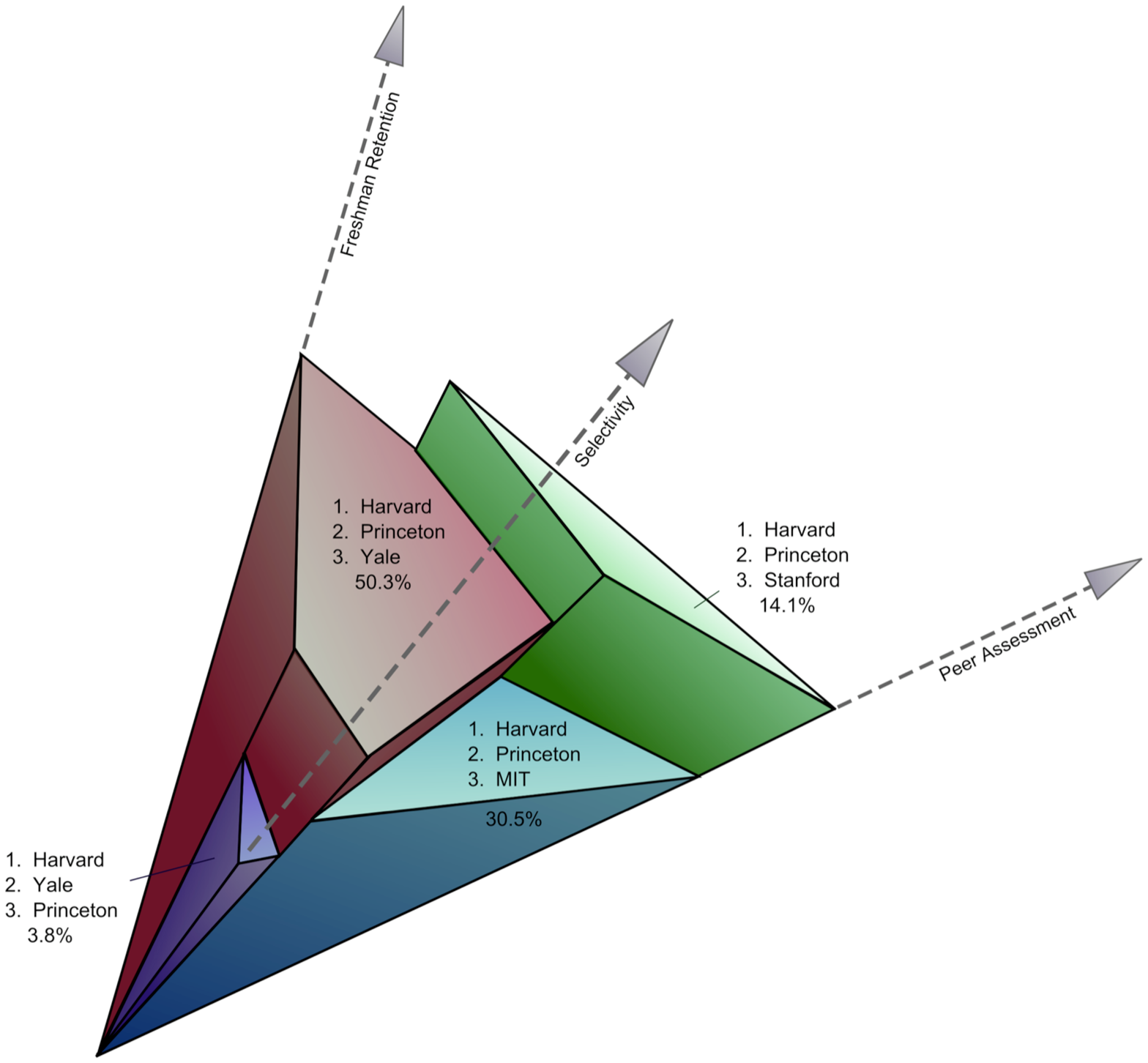}
\vskip -1in
\caption{University rankings determined by choice of weights. The top three universities as determined by the choice of weights for three categories. The four cones comprise 99\% of the total volume. }
\end{figure}

  We then sampled functionals $c$ from the non-negative first orthant using the Gaussian density, and computed the rankings that the functionals induced.  
  Geometrically, this gives approximate volumes of the normal cones of the $k$-ordertopes.  $2,000,000$ functionals were sampled and we computed various summaries of the induced rankings.  We recorded how many times each university was assigned each rank between 1 and 130, and how many times university $i$ was ranked higher
than university $j$ for all pairs $(i,j)$.  

  For all 130 universities we computed ranking intervals, which are given below:
  \begin{tabular}{lrc}
   University & USNWR Ranking & Ranking CI \\
   \hline \\
Princeton University (NJ)    & 1 &      [1,4] \\
Harvard University (MA)    & 2 &        [1,4] \\
Yale University (CT)    & 3 &   [1,5] \\
Stanford University (CA)    & 4 &       [3,8] \\
University of Pennsylvania    & 5 &     [3,10] \\
California Institute of Technology    & 5 &     [2,21] \\
Massachusetts Institute of Technology    & 7 &  [2,15] \\
Duke University (NC)    & 8 &   [2,11] \\
Columbia University (NY)    & 9 &       [6,13] \\
University of Chicago    & 9 &  [7,18] \\
Dartmouth College (NH)    & 11 &        [3,12] \\
Washington University in St Louis    & 12 &     [5,19] \\
Cornell University (NY)    & 12 &       [9,15] \\
Brown University (RI)    & 14 &         [9,18] \\
Northwestern University (IL)    & 14 &  [11,19] \\
Johns Hopkins University (MD)    & 14 &         [4,18] \\
Rice University (TX)    & 17 &  [13,21] \\
Emory University (GA)    & 17 &         [14,35] \\
Vanderbilt University (TN)    & 19 &    [16,26] \\
University of Notre Dame (IN)    & 19 &         [6,22] \\
University of California Berkeley *    & 21 &    [14,36] \\
Carnegie Mellon University (PA)    & 22 &       [19,36] \\
University of Virginia *    & 23 &      [20,33] \\
Georgetown University (DC)    & 23 &    [15,26] \\
University of California Los Angeles *    & 25 &         [17,35] \\
University of MichiganAnn Arbor *    & 25 &     [18,37] \\
University of Southern California    & 27 &     [20,35] \\
University of North Carolina Chapel Hill *    & 28 &     [26,35] \\
Tufts University (MA)    & 28 &         [22,33] \\
Wake Forest University (NC)    & 30 &   [12,37] \\
Lehigh University (PA)    & 31 &        [17,38] \\
Brandeis University (MA)    & 31 &      [23,37] \\
College of William and Mary (VA) *    & 33 &    [27,48] \\
New York University    & 34 &   [32,78] \\
University of Rochester (NY)    & 35 &  [22,59] \\
Georgia Institute of Technology *    & 35 &     [25,73] \\
Boston College    & 35 &        [28,43] \\
University of Wisconsin Madison *    & 38 &      [31,58] \\
University of California San Diego *    & 38 &   [31,81] \\
University of Illinois Urbana - Champaign *    & 38 &    [31,52] \\
Case Western Reserve University (OH)    & 41 &  [26,101] \\
University of Washington *    & 42 &    [30,53] \\
University of California Davis *    & 42 &       [37,72] \\

\end{tabular}

\begin{tabular}{lrc}
   University & USNWR Ranking & Ranking CI \\
   \hline \\
Rensselaer Polytechnic Institute (NY)    & 44 &         [38,58] \\
University of Texas Austin *    & 44 &   [36,71] \\
University of California Santa Barbara *    & 44 &       [36,64] \\
University of California Irvine *    & 44 &      [41,88] \\
Pennsylvania State University University Park *    & 48 &        [1,59] \\
University of Florida *    & 49 &       [42,66] \\
Syracuse University (NY)    & 50 &      [31,60] \\
Tulane University (LA)    & 50 &        [35,109] \\
Yeshiva University (NY)    & 52 &       [26,75] \\
University of Miami (FL)    & 52 &      [38,74] \\
Pepperdine University (CA)    & 54 &    [40,86] \\
George Washington University (DC)    & 54 &     [49,101] \\
University of Maryland College Park *    & 54 &  [47,81] \\
Ohio State University Columbus *    & 57 &       [42,73] \\
Boston University    & 57 &     [52,108] \\
Rutgers, Univ. of New Jersey New Brunswick *    & 59 &       [43,73] \\
University of Pittsburgh *    & 59 &    [43,74] \\
University of Georgia *    & 59 &       [56,83] \\
Texas A and M University College Station *    & 62 &         [52,85] \\
Worcester Polytechnic Institute (MA)    & 62 &  [45,96] \\
University of Connecticut *    & 64 &   [38,72] \\
Purdue UniversityWest Lafayette (IN) *    & 64 &        [43,95] \\
University of Iowa *    & 64 &  [55,89] \\
Fordham University (NY)    & 67 &       [30,83] \\
Miami University Oxford (OH) *    & 67 &         [33,98] \\
Clemson University (SC) *    & 67 &     [42,75] \\
Southern Methodist University (TX)    & 67 &    [58,86] \\
Univ. of Minnesota Twin Cities *    & 71 &  [56,120] \\
Virginia Tech *    & 71 &       [40,89] \\
University of Delaware *    & 71 &      [46,75] \\
Michigan State University *    & 71 &   [33,95] \\
Stevens Institute of Technology (NJ)    & 75 &  [43,115] \\
Baylor University (TX)    & 75 &        [30,94] \\
Colorado School of Mines & 75 &         [46,120] \\
Indiana University Bloomington *    & 75 &       [41,115] \\
Brigham Young University Provo (UT)    & 79 &    [58,113] \\
University of California Santa Cruz *    & 79 &  [73,127] \\
University of Colorado Boulder *    & 79 &       [77,123] \\
St Louis University    & 82 &   [72,91] \\
SUNY Binghamton *    & 82 &      [49,127] \\
Marquette University (WI)    & 82 &     [49,100] \\
SUNY College of Environmental     & 85 &  [35,105] \\
Science and Forestry *    &  &    \\
North Carolina State Univ. Raleigh *    & 85 &      [42,105] \\

\end{tabular}

\begin{tabular}{lrc}
   University & USNWR Ranking & Ranking CI \\
   \hline \\

University of Denver    & 85 &  [68,107] \\
American University (DC)    & 85 &      [79,111] \\
Iowa State University *    & 85 &       [70,99] \\
University of Kansas *    & 85 &        [77,122] \\
University of Alabama *    & 91 &       [46,113] \\
University of Missouri Columbia *    & 91 &      [93,126] \\
University of Nebraska Lincoln *    & 91 &       [75,105] \\
University of Tulsa (OK)    & 91 &      [58,130] \\
Clark University (MA)    & 91 &         [59,105] \\
Auburn University (AL) *    & 96 &      [87,111] \\
SUNY Stony Brook *    & 96 &     [87,125] \\
University of Tennessee *    & 96 &     [78,128] \\
University of Vermont *    & 96 &       [75,109] \\
University of Arizona *    & 96 &       [84,129] \\
University of the Pacific (CA)    & 96 &        [64,118] \\
University of California Riverside *    & 96 &   [88,124] \\
Howard University (DC)    & 96 &        [33,111] \\
Illinois Institute of Technology    & 96 &      [79,130] \\
Northeastern University (MA)    & 96 &  [90,122] \\
University of Massachusetts Amherst *    & 96 &  [78,117] \\
University of San Diego    & 107 &      [89,118] \\
University of New Hampshire *    & 108 &        [56,126] \\
University of Oklahoma *    & 108 &     [94,124] \\
Drexel University (PA)    & 108 &       [79,125] \\
Texas Christian University    & 108 &   [77,124] \\
Ohio University *    & 112 &    [40,130] \\
University of Dayton (OH)    & 112 &    [56,122] \\
University of Oregon *    & 112 &       [78,129] \\
University of South Carolina Columbia *    & 112 &       [88,122] \\
Florida State University *    & 112 &   [97,130] \\
Loyola University Chicago    & 112 &    [87,129] \\
University of Missouri Rolla *    & 118 &        [96,130] \\
University at Buffalo SUNY *    & 118 &  [96,128] \\
Washington State University *    & 118 &        [86,124] \\
Samford University (AL) & 118 &         [97,129] \\
Catholic University of America (DC)    & 122 &  [111,130] \\
University of Kentucky *    & 122 &     [91,129] \\
New Jersey Institute of Technology *    & 124 &         [88,130] \\
Clarkson University (NY)    & 124 &     [87,128] \\
Arizona State University *    & 124 &   [114,130] \\
University of Arkansas *    & 124 &     [103,130] \\
Colorado State University *    & 124 &  [104,130] \\
Kansas State University *    & 124 &    [94,128] \\
Michigan Technological University *    & 124 &  [108,129] \\

\end{tabular}

  \subsection{Validation of results on random university data}  

  Since we could not perfectly reverse-engineer the university score vectors from the USNWR data, we simulated data to verify that our results are
  typical for rankings of seven-dimensional data.  To generate the $i$th university's data $u_i$, we first sampled an intrinsic goodness $t_i \in {\bf R}$ from the normal distribution.  We then set

  \[
     u_{ij} = p t_i + (1 - p) \epsilon_{ij}
  \]
  where each $\epsilon_{ij}$ is sampled from the normal distribution and $p \in [0,1]$ is a fixed constant.  Thus the seven category scores in $u_i$ are correlated and influenced by the intrinsic goodness $t_i$, but also fluctuate somewhat from one category to the next, due to the added $\epsilon_{ij}$.  If $p = 1$ then there is no
  fluctuation between categories and there is only one possible ranking of the universities for non-negative choice of weights $c \neq 0$. If $p = 0$ then the cateogires are completely uncorrelated.  We chose $p = 0.6$, because this gave higher correlation between category rankings on average, compared to the correlations observed in USNWR category rankings.  

  Despite our conservative choice of $p = 0.6$, we still observed considerable variability in the possible overall rankings. 
  Across 50 trials, the average ranking interval width for the 130 universities ranged between 20.02 and 28.19, and was 22.90 on average.  
  
  This shows that volatility in rankings is typical for seven dimensional data, even when the seven coordinates are 
  influenced by a common variable.  Intuitively, this problem becomes worse as the dimension of the data increases.
.

\pagebreak

%%%%%%%%%%%%%%%%%%%%%%%%%%%%%%%%%%%%%%%%%%%%%%%%%%%%%%%%%%%%%%%%%%%%%%%%%

\end{document}